\documentclass[12pt]{amsart}
\usepackage{amsmath,amssymb}
\usepackage{amsfonts}
\usepackage{amsthm}
\usepackage{latexsym}
\usepackage{graphicx}

\def\ddbar{\partial\bar\partial}

\def\iddbar{\sqrt{-1}\ddbar}
\def\ii{\sqrt{-1}}
\def\rank{\mbox{rank}}

\def\C{\mathbb{C}}

\def\vv<#1>{\langle#1\rangle}

\def\XXint#1#2{\setbox0=\hbox{$#1{#2}{\int}$}{#2}\kern-.5\wd0 }

\def\XXint#1#2#3{{\setbox0=\hbox{$#1{#2#3}{\int}$}
     \vcenter{\hbox{$#2#3$}}\kern-.5\wd0}}

\def\vv<#1>{\langle#1\rangle}

\def\ii{\sqrt{-1}}

\def\be{\begin{equation}}
\def\ee{\end{equation}}

\def\C{\Bbb C}

\def\C{\Bbb C}
 \def\ii{\sqrt{-1}}

\newtheorem{thm}{Theorem}[section]

\newtheorem{cor}{Corollary}[section]
\newtheorem{conjecture}{Conjecture}
\theoremstyle{definition}

\theoremstyle{remark}

\newtheorem{rem}{Remark}[section]

\numberwithin{equation}{section}

\begin{document}
\title{A note on Wu-Zheng's splitting conjecture}

\keywords{ K\"ahler manifolds, bisectional curvature, splitting theorem}

\begin{abstract} Cao's splitting theorem says that for any
complete K\"ahler-Ricci flow $(M,g(t))$ with $t\in [0,T)$, $M$
simply connected and  nonnegative bounded holomorphic bisectional
curvature, $(M,g(t))$ is holomorphically  isometric to $\C^k\times
(N,h(t))$ where $(N,h(t))$ is a K\"ahler-Ricci flow with positive
Ricci curvature for $t>0$. In this article, we show that $k=n-r$
where $r$ is the Ricci rank of the initial metric. As a corollary,
we also confirm a splitting conjecture of Wu-Zheng when curvature is
assumed to be bounded.
\end{abstract}

\renewcommand{\subjclassname}{\textup{2010} Mathematics Subject Classification}
 \subjclass[2010]{Primary 53C44 ; Secondary  53C55}
\author{Chengjie Yu$^1$}

\address{Department of Mathematics, Shantou University, Shantou, Guangdong, China }
\email{cjyu@stu.edu.cn}

\thanks{$^1$Research partially supported by the National Natural Science Foundation of
China (11001161) and (10901072).}
\date{Nov. 2010}
\maketitle \markboth{Chengjie Yu}
 {Wu-Zheng's splitting conjecture}
 \section{Introduction}
 In \cite{cao}, Cao proved the following splitting theorem.
\begin{thm}[Cao H.-D.]\label{thm-cao}
Let $g(t)$ with $t\in[0,T)$ be a complete solution to the K\"ahler-Ricci flow
\begin{equation}\label{eqn-K-R}
\frac{d}{dt}g_{i\bar j}=-R_{i\bar j}
\end{equation}
on a complex manifold $M^n$ and the initial metric $g(0)$ has nonnegative and bounded holomorphic bisectional curvature.
Then, $\rank\ R_{i\bar j}(x,t)$ is independent of $x\in M$ and $t\in (0,T)$. Moreover, $\ker R_{i\bar j}(x,t)$ is
independent $t>0$ and is a parallel distribution for all $t>0$. In particular, if $M$ is simply connected,
$(M,g(t))$ is holomorphically isometric to $\C^{n-k}\times (N^{k},h(t))$ where $h(t)$ with $t\in [0,T)$ is a complete solution to the
K\"ahler-Ricci flow \eqref{eqn-K-R} on $N$ with positive Ricci curvature for $t>0$.
\end{thm}

Let $(M,g)$ be a complete K\"ahler manifold and $r=\max_{x\in
M}\rank\ R_{i\bar j}(x)$. Then $r$ is called the Ricci rank of the
K\"ahler metric $g$. In this article, we shown that the number $k$
in Cao's splitting theorem equals the Ricci rank of the initial
metric $g(0)$.
\begin{thm}\label{thm-Ricci-rank}
Let settings be the same as in Theorem \ref{thm-cao} and $r$ be the Ricci rank of the initial metric $g(0)$. Then
\begin{equation}
\rank\ R_{i\bar j}(x,t)=r
\end{equation}
for any $x\in M$ and $t\in (0,T)$.
\end{thm}

In \cite{Wu-Zheng}, Wu-Zheng presented the following splitting conjecture:
\begin{conjecture}[Wu-Zheng]
Let $M^n$ be a complete K\"ahler manifold with nonnegative holomorphic bisectional curvature. Then, its universal covering manifold $\tilde M$ is
holomorphically isometric to $\C^{n-r}\times N^r$, where $r$ is the Ricci rank.
\end{conjecture}
This is a noncompact generalization of the splitting theorem of Howard-Smyth-Wu \cite{HSW} on compact K\"ahler manifolds with
nonnegative holomorphic bisectional curvature which is also related to a conjecture proposed by Yau (Ref. \cite{Wu-Zheng,Wu-Zheng-2}). In \cite{Wu-Zheng} and \cite{Wu-Zheng-2}, Wu-Zheng showed that the conjecture is true for complete
K\"ahler manifolds with nonnegative holomorphic bisectional curvature, real analytic K\"ahler metric and Ricci rank $\leq 2$. The technique used by Zheng-Wu is delicate. They performed some careful analysis on the distribution defined by the kernel of the Ricci tensor. Because the distribution is not defined all over the manifold, they have to require that the metric is real analytic, so that the distribution can be paralleled extended to the whole manifold. As a direct corollary of Theorem \ref{thm-cao} and Theorem \ref{thm-Ricci-rank}, we can
see that the conjecture is true if curvature is assumed to be bounded.
\begin{thm}\label{thm-wu-zheng-splitting}
Let $M^n$ be a complete K\"ahler manifold with nonnegative and bounded holomorphic bisectional curvature. Then, its universal covering manifold $\tilde M$ is
holomorphically isometric to $\C^{n-r}\times N^r$, where $r$ is the Ricci rank.
\end{thm}

In \cite{Ni-Tam}, Ni-Tam got splitting results of complete K\"ahler manifolds with nonnegative holomorphic bisectional curvature in a
similar spirit of Cao's splitting theorem. However, Ni-Tam's splitting results are with respect to function theory on the
K\"ahler manifold. For example, they got the following result.
\begin{thm}[Ni-Tam]\label{thm-ni-tam}
Let $M$ be a complete simply connected noncompact K\"ahler manifold with nonnegative holomorphic bisectional curvature and supporting a
smooth strictly plurisubharmonic function $u$ on $M$ with bounded gradient. Then $M=\C^l\times M_1\times M_2$ isometrically and
holomorphically for some $l\geq 0$, where $M_1$ and $M_2$ are complete noncompact K\"ahler manifolds with nonnegative
holomorphic bisectional curvature such that any polynomial growth holomorphic function on $M$ is independent of the factor
$M_2$, and any linear growth holomorphic function is independent of the factors $M_1$ and $M_2$.
\end{thm}
As a direct corollary of Theorem \ref{thm-wu-zheng-splitting}, we know that $l$ equals the Ricci rank of the manifold if
we further assume that the curvature of the manifold is bounded.
\begin{cor}
Let settings be the same as in Theorem \ref{thm-ni-tam} and further assume that the curvature of $M$ is bounded. Then
$l$ equals to the Ricci rank of the manifold in the conclusion of Theorem \ref{thm-ni-tam}.
\end{cor}
\section{proof of Theorem \ref{thm-Ricci-rank}}
In this section, using a result of Ferus \cite{Ferus} on the completeness of the leaves of the nullity foliation of the curvature
operator, we give a proof of Theorem \ref{thm-Ricci-rank}.
\begin{proof}[Proof of Theorem \ref{thm-Ricci-rank}]
Let $\Omega=\{x\in M\ |r=\rank\ R_{i\bar j}(x,0)\}$ and $L_x=\ker R_{i\bar j}(x,0)$ for each $x\in \Omega$. It is shown in
\cite{Wu-Zheng} that $L$ is a totally geodesic and flat foliation on $\Omega$. In \cite{Ferus}, Ferus showed that each leaf of the foliation
$L$ is complete. Let $F_x=\ker R_{i\bar j}(x,t)$ for each $x\in M$ and $t>0$. It is shown by Cao's splitting theorem (Theorem \ref{thm-cao}) that
$F$ is independent of $t$ and parallel for each $t>0$. By continuity(Ref. \cite{Shi}), we know that $F\subset L$ on $\Omega$ and $F$ is
also parallel with respect to $g(0)$. It suffices to show that $L\subset F$ on $\Omega$.

Since $L$ is totally geodesic, flat and with complete leaves, we
know that a leaf of $L$ is a totally geodesic and isometric
holomorphic immersion $\phi:\C^{n-r}\to M$. Note that on a
K\"ahler-Ricci flow,
\begin{equation}
\frac{d}{dt}Ric=\iddbar R
\end{equation}
where $Ric=\ii R_{i\bar j}dz\wedge dz^{\bar j}$ is the Ricci form. Therefore, for each $t>0$,
\begin{equation}
Ric(t)=Ric(0)+\iddbar\int_0^{t}R(s)ds
\end{equation}
and
\begin{equation}
\phi^*Ric(t)=\phi^*Ric(0)+\iddbar\int_0^{t}R(s)(\phi)ds.
\end{equation}
Note that $\phi^*Ric(0)=0$ since $\phi:\C^{n-r}\to M$ is a leaf of $L$. Hence
\begin{equation}\label{eqn-Ric-t}
\phi^*Ric(t)=\iddbar\int_0^{t}R(s)(\phi)ds.
\end{equation}
Note that $Ric(t)$ is nonnegative, so $\int_0^{t}R(s)(\phi)ds$ is a
bounded pluri-subharmonic function on $\C^{n-r}$. By Liouville
property of pluri-subharmonic functions on complex Euclidean spaces
(See \cite[Theorem 0.2]{Ni-Tam}), we know that
$\int_0^{t}R(s)(\phi)ds$ is a constant function. Substituting this
back to \eqref{eqn-Ric-t}, we know that $\phi^*Ric(t)=0$ for $t>0$.
Therefore $L\subset F$ on $\Omega$.
\end{proof}
\begin{rem}
Instead of using paralleled translation to extend the Ricci kernel, here we use the K\"ahler-Ricci flow to extend the Ricci kernel. 
\end{rem}

\end{document}